\documentclass[journal,onecolumn]{IEEEtran}
\usepackage{cite}

\ifCLASSINFOpdf
\else
\fi

\usepackage{graphics} 
\usepackage{epsfig} 
\usepackage{mathptmx} 
\usepackage{times} 
\usepackage{amsmath} 
\usepackage{amssymb}  
\usepackage{bbm}
\usepackage{chemfig,siunitx}
\usepackage{textcomp}
\usepackage{subfigure}

\newtheorem{theorem}{Theorem}

\newtheorem{lemma}{Lemma}
\newtheorem{rem}{Remark}
\makeatletter
\newcommand{\onetagright}{\tagsleft@false}
\makeatother

\begin{document}
%
\title{Input Novelty as a Control Metric for Time Varying Linear Systems}
%
%
%

\author{Gautam~Kumar,~\IEEEmembership{Member,~IEEE,}
        Delsin~Menolascino,~\IEEEmembership{Student Member,~IEEE,}
        and~ShiNung~Ching,~\IEEEmembership{Member,~IEEE}
\thanks{G. Kumar, D. Menolascino, and S. Ching are with the Department
of Electrical and Systems Engineering, Washington University in St. Louis, St. Louis, MO, 63130 USA e-mail: (shinung@ese.wustl.edu).}}
\maketitle

\begin{abstract}
This paper introduces a framework for quantitative characterization of the controllability of time-varying linear systems (or networks) in terms of input novelty.
The motivation for such an approach comes from the study of biophysical sensory networks in the brain, wherein responsiveness to both energy and salience (or novelty) are presumably critical for mediating behavior and function.  Here, we use an inner product to define the angular separation of the current input with respect to the past input history.  Then, by constraining input energy, we define a non-convex optimal control problem to obtain the minimally novel input that effects a given state transfer.  We provide analytical conditions for existence and uniqueness in continuous-time, as well as an explicit closed-form expression for the solution.  In discrete time, we show that a relaxed convex optimization formulation provides the global optimal solution of the original non-convex problem.  Finally, we show how the minimum novelty control can be used as a metric to study control properties of large scale recurrent neuronal networks and other complex linear systems.  In particular, we highlight unique aspects of a system's controllability that are captured through the novelty-based metric.  The result suggests that a multifaceted approach, combining energy-based analysis with other specific metrics, may be useful for obtaining more complete controllability characterizations.

\end{abstract}


%
\IEEEpeerreviewmaketitle

\section{Introduction}
In its most basic form, the systems-theoretic notion of controllability carries a binary definition: a dynamical system either is, or is not, controllable, with respect to its exogenous inputs. Naturally, such a notion has the deficiency of not grading the ease or difficulty associated with particular control objectives. To obviate this issue, consistent research effort has been directed at the characterization of controllability using systems-theoretic metrics. Roughly, these metrics can be grouped into two categories

\begin{enumerate}
	\item Those that characterize the minimum energy parametric perturbations -- termed controllability radii -- that result in a loss of controllability \cite{Hu2001,HD04}.  These are related to basic characterizations of the robustness of linear systems \cite{hinrichsen1989real}.
	\item Those that characterize the controllability of a system in terms of the minimum energy excitation required to effect a desired state transfer \cite{YRCHL12,PZB14,FZ14,CSL14}.  
\end{enumerate}
The latter, in particular, is a natural paradigm that is directly related to the Kalman rank condition and, specifically, the controllability gramian, used to ascertain the controllability of linear systems \cite{KS72}.  Recently, energy-based controllability metrics have been successfully used in the emerging domain of network science to assess the putative controllability of large-scale networks with linear node and connection dynamics \cite{FZ14,PZB14}. However, for complex networks in general and for biological neuronal networks in particular, there are certain questions related to the overall system dynamics that are not fully captured through an energy-based metric alone.

\begin{figure}
	\centering
		\includegraphics[width=0.40\textwidth]{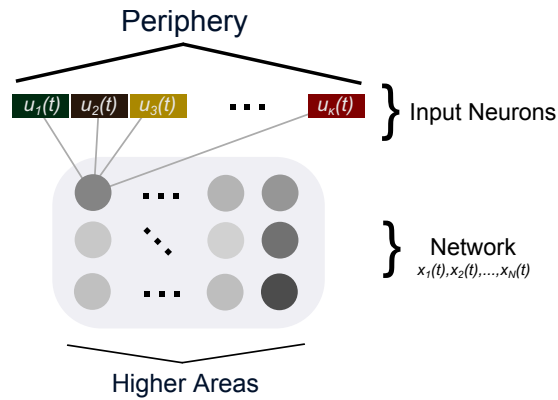}
	\caption{Prototypical structure of a sensory neuronal network. Sensory neurons are tuned to features from the sensory periphery. These neurons project excitation onto a network that performs intermediate transformations on the afferent excitation en route to higher brain regions.}
	\label{fig1}
\end{figure}

We appeal, specifically, to the domain of neural coding and the dynamics of sensory networks in the brain. Consider the simple, prototypical `feedforward-type' layered model of a sensory network shown in Figure \ref{fig1}, wherein sensory neurons are tuned to a high dimensional feature space (i.e., environmental variables from the sensory periphery; say, different molecules corresponding to tastes). Those sensory neurons impinge on a complex, interconnected, recurrent sensory network that performs intermediate transformations on the incoming activity, en route to higher brain areas.

One may put forth a supposition that the controllability of such a sensory network, with respect to the afferent input from the sensory neurons, is critical for facilitating perception and behavior.  For instance, a network that requires an enormous excitation in order to effect any state trajectory would ostensibly be quite poor in mediating perception of small changes in the environment. 

But as much as energy is important in such mediation, \textit{orientation}, i.e., the alignment of an input with certain features, and \textit{novelty}, the difference in orientation of an input from past inputs, may be equally so. Indeed, a weak, but highly \textit{novel} input may be more easily perceived than an intense, but more familiar, stimulus \cite{Downar2002}. The ability to assess the responsiveness of neuronal networks to novelty -- at a particular moment in time, relative to past inputs -- has immediate implications in the analysis and control of physiological neuronal network dynamics in different behavioral and clinical regimes \cite{Ching2013,Ching2012,Lepage2013a}.

\begin{figure}[t]
\centering%
\includegraphics[width=0.48\textwidth]{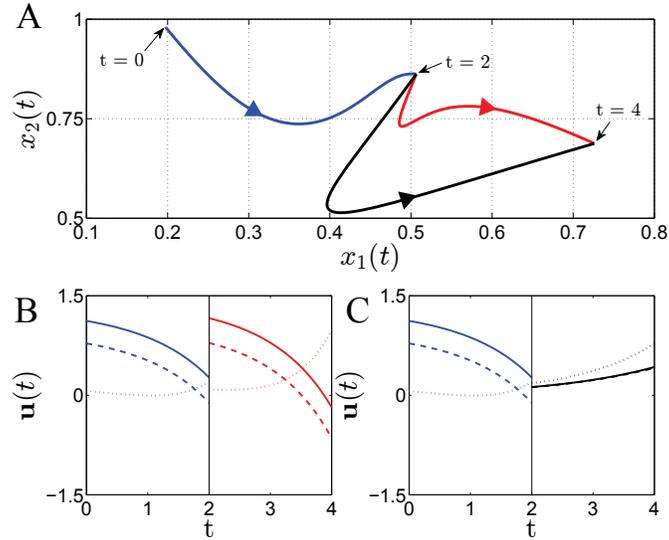}
\caption{Minimum novelty control vs. minimum energy control: (A) The trajectory (blue) brings the system from an initial state on intermediate state at $t=2s$.  Subsequently, two trajectories are contrasted in the phase-plane for the minimum novelty control (red) and the minimum energy control (black). (B) The minimally novel inputs (from $t=2s$ to $t=4s$) (red) designed using our approach in this paper. (C) The inputs corresponding to the minimum energy trajectory (from $t=2s$ to $t=4s$), (black).}
\label{fig2}
\end{figure}

In this paper, we seek a graded quantification of controllability (or, if one prefers, reachability) of linear time-varying systems (networks) in terms of input novelty. In particular, we ask how responsive are the state (node) trajectories to inputs that differ in orientation from those that have previously been applied. Figure \ref{fig2} illustrates this basic notion for a simple two-dimensional linear system with a three-dimensional input. A particular input drives the system from a point in the phase space at $t=0$ to an intermediate point at $t=2$; from this point emerge two trajectories, both of which reach a common endpoint; one minimizes input novelty (note the similarity between the input over $t\in [0,2]$ and that over $t\in [2,4]$), while the other minimizes energy.  

We proceed to formalize this notion and develop its utility in control analysis.  Our contributions are as follows:
\begin{enumerate}
\item  We define the notion of input novelty through the use of an inner product evaluated in the input feature-space
\item  We analytically derive, for continuous, linear time-varying systems, the \textit{minimum novelty control} that effects a desired state transfer by formulating a non-convex optimization problem.  The problem seeks the minimum angular separation, defined in terms of the inner product, required in order to create a desired change in the network trajectory, constrained by a fixed average input energy.
\item  For discrete-time systems, we obtain the minimum novelty control through a provably exact convex relaxation of the original non-convex problem. The relaxed problem can be solved efficiently using well established convex programming solvers.
\item We use the minimally novel control to construct a metric which can be used to characterize large-scale linear systems or networks.  Importantly, we demonstrate that this novelty-based metric provides a distinct characterization of controllability compared to energy-based control metrics.  The result suggests that a broader, multifaceted set of metrics may be useful for completely characterizing the controllability of such systems.
\end{enumerate}


The remainder of the paper is organized as follows. In Section \ref{sec2}, we describe the mathematical notation used in this paper. In Section \ref{sec3}, we introduce the inner product-based input novelty measure for continuous-time, linear time-varying systems and formulate a non-convex optimal control problem that minimizes this novelty under the constraint of a fixed average input energy.  We establish the existence and the uniqueness of a global optimal solution of the control problem and derive a closed-form expression for the minimally novel input.
In Section \ref{sec4}, we derive analogous results for discrete time linear time-varying systems with the additional result of exact relaxation of the non-convex novelty minimization problem.
Finally, in Section \ref{sec5}, we demonstrate the use of the minimally novel input as a control metric to systematically characterize large scale linear networks. The paper concludes with a summary and discussion of future work.

\section{Mathematical Notation}\label{sec2}
\noindent Most notation is standard and will be introduced as the results are developed. We use lower-case letters to represent scalars, boldface lower-case letters to represent vectors, capital letters to represent matrices. Exceptions are $T$, $T^{*}$, $\mathbb{J}(T)$, $\mathbb{J}_{1}(T)$ and $\mathbb{J}(p)$, which we represent as scalars. We use $\mathbb{R}^{n\times 1}$ to denote the space of $n$- dimensional vectors with their elements as real numbers. Similarly, we use $\mathbb{R}^{n\times m}$ and $\mathbb{R}_{+}^{n\times m}$ to denote the space of $n\times m$ dimensional matrices with real and non-negative real entries, respectively. $\|\mathbf{x}\|_{2}$ is the Euclidean norm of the vector $\mathbf{x}$. $\mathbf{x}^{'}$ is the transpose of a vector $\mathbf{x}$. We use $\mathbf{x}(t)$ and $\mathbf{x}(k)$ to represent the value of $\mathbf{x}$ at continuous time $t$ and discrete time $k$ respectively. $\mathbf{x}_{r}$, $\mathbf{x}_{0}$ and $\mathbf{x}_{f}$ are constant $n$-dimensional vectors. $A^{-1}$ is the inverse of a matrix $A$ and $(A)_{ij}$ represents the $(ij)^{th}$ element of $A$.

\section{Continuous-Time, Linear Dynamical Networks}\label{sec3}

\subsection{Input Novelty based Controllability Metric}
\noindent We consider a linear, time-varying system with dynamics of the form
\begin{equation}
	\frac{\rm{d}\mathbf{x}(t)}{\rm{d}t} = A(t) \mathbf{x}(t) + B(t) \mathbf{u}(t)
\label{eq1}
\end{equation} 
Here $\mathbf{x}(t)\in \mathbb{R}^{n\times 1}$ represents the state of the system at time $t$, $A(t)\in \mathbb{R}^{n\times n}$ describes the time-varying dynamics, $B(t)\in \mathbb{R}^{n\times m}$ is the input matrix, and $\mathbf{u}(t)\in \mathbb{R}^{m\times 1}$ is the input to the system.  Without loss of generality, we say that (\ref{eq1}) describes the time evolution of (time-varying) linear networks in the presence of external inputs.  

We consider an input $\mathbf{v}(t-T) \in \mathbb{R}^{m\times 1}$, $t\in[0,T]$, with total energy $T\gamma_{v}$, i.e.
\begin{equation}
\frac{1}{T}\int_{0}^{T}\|\mathbf{v}(t-T)\|_{2}^{2}\rm{d}t = \gamma_{v} 
\label{eq2}
\end{equation}
Further, we assume that $\mathbf{v}(t-T)$ drives $\mathbf{x}(t)$ from $\mathbf{x}_{r}$ to $\mathbf{x}_{0}$ subject to the dynamics (\ref{eq1}). Here, $T>0$ is a constant and $\gamma_{v}>0$ is the average per-time input energy. We introduce the inner product-based measure   

\begin{equation}
	\mathbb{J}(T) = \frac{1}{T\sqrt{\gamma_{v}\gamma_{u}}}\int_{0}^{T} \mathbf{v}^{'}(t-T)\mathbf{u}(t)\rm{d}t
\label{eq3}
\end{equation}

\noindent where 
\begin{equation}
\frac{1}{T}\int_{0}^{T}\|\mathbf{u}(t)\|_{2}^{2}\rm{d}t = \gamma_{u},
\label{eq4}
\end{equation}
\noindent to quantify the novelty of a subsequent input $u(t)$, $t\in[0,T]$ (which drives $\mathbf{x}(t)$ from $\mathbf{x}_{0}$ to $\mathbf{x}_{f}$), relative to $v(t-T)$.
For a fixed input energy, $\mathbb{J}(T)$ measures the change in orientation (thus novelty) of two consecutive inputs trajectories of equal lengths (in time). 

\begin{rem}\label{rem1}
It is readily evident that $\mathbb{J}(T)\in [-1,1]$. 
\end{rem}

\begin{rem}\label{rem2}
It follows from (\ref{eq3}) that the novelty of $\mathbf{u}(t)$ relative to $\mathbf{v}(t-T)$ decreases as $\mathbb{J}(T)$ increases and is minimum when $\mathbb{J}(T) = 1$ i.e. when $\mathbf{u}(t) = \alpha \mathbf{v}(t-T)$, $\alpha\in \mathbb{R}$, for all $t\in[0,T]$.
\end{rem}

\begin{rem}\label{rem3}
	 We observe that, due to the energy normalization in (\ref{eq2}) and (\ref{eq4}),
\begin{equation}
\frac{1}{T}\int_{0}^{T} \|\mathbf{v}(t-T)-\mathbf{u}(t)\|_{2}^{2}\rm{d}t =  (\gamma_{v}+\gamma_{u}-2\sqrt{\gamma_{v}\gamma_{u}} \mathbb{J}(T))
\label{eq5}
\end{equation}
\noindent Thus, the average Euclidean distance, i.e. the left hand side of (\ref{eq5}), between two inputs can be used to generate the input novelty in our context.   
\end{rem}

\subsection{Minimum Novelty Problem}
\noindent From the conceptual formulation introduced above, we can develop a control problem to design the minimally novel input $\mathbf{u}(t), t\in[0,T]$ such that a desired change in the state of the system can be achieved under the constraint of fixed energy subject to (\ref{eq1}). Specifically, we formulate the following optimal control problem:

\begin{subequations}
\begin{align}
    \min_{\substack{\mathbf{u}(t) \\ t\in[0,T]}}
        & \quad  -\mathbb{J}(T)  \label{eq6a}\\
    \textrm{s.t.} 
        & \qquad \frac{1}{T} \int_{0}^{T}\| \mathbf{u}(t)\|_{2}^{2}\rm{d}t = \gamma_{u} \label{eq6b}\\
        & \quad  \mathbf{x}_{f} = \Phi(T,0)\mathbf{x}_{0} + \int_{0}^{T}\Phi(T,t)B(t)\mathbf{u}(t)\rm{d}t \label{eq6c} 
\end{align}
\label{eq6}
\end{subequations}

\noindent Here, the state-transition matrix $\Phi(t,\tau)$ is given by
\begin{equation}
	\frac{\partial \Phi(t,\tau)}{\partial t} = A(t)\Phi(t,\tau)
\label{eq7}
\end{equation} 
\noindent with $\Phi(\tau,\tau) = I$, where $I$ is the $n\times n$ identity matrix.   

It should be noted here that the constraint (\ref{eq6c}) is obtained by integrating (\ref{eq1}) with respect to $t$ over the period of $[0,T]$. Immediately, we note that the quadratic equality constraint (\ref{eq6b}) makes the optimization problem (\ref{eq6}) non-convex. Furthermore, we note that our optimal control problem formulation (\ref{eq6}) is different from the classical minimum effort problems \cite{K04,NQN13} where the $L^{1}$-norm of control inputs is minimized under the constraints of explicit lower and upper bounds on the inputs. 

\subsection{Analytical Results}
\noindent We derive conditions for the existence of a unique global optimal solution of the non-convex optimization problem (\ref{eq6}). Based on this, we provide a closed-form expression for the optimal $\mathbf{u}(t), t\in[0,T]$.

\subsubsection{Existence of a Minimally Novel Input}

\begin{theorem}\label{thm1}
	A solution of the non-convex optimization problem (\ref{eq6}) exists if 
\begin{equation}
	T > \max{\{\frac{\mathbf{s}^{'}(T)W_{c}^{-1}(T)\mathbf{s}(T)}{\gamma_{v}},\frac{\mathbf{r}^{'}(T)W_{c}^{-1}(T)\mathbf{r}(T)}{\gamma_{u}}\}}
\label{eq8}
\end{equation}

\noindent where 
\begin{subequations}
\begin{equation}
\mathbf{s}(T) = \int_{0}^{T}\Phi(T,t)B(t)\mathbf{v}(t-T)\rm{d}t
\label{eq9a}
\end{equation}
\begin{equation}
\mathbf{r}(T) = \mathbf{x}_{f}-\Phi(T,0)\mathbf{x}_{0}
\label{eq9b}
\end{equation}
\label{eq9}
\end{subequations}

\noindent Here, $W_{c}(T)$\footnote{For notational simplicity, dependence on the initial time (here, $t_0 = 0$) is implicit.} is the usual controllability gramian at time $T$ \cite{KS72} and is defined as 
\begin{equation}
	W_{c}(T) = \int_{0}^{T}\Phi(T,t)B(t)B^{'}(t)\Phi^{'}(T,t)\rm{d}t
\label{eq10}
\end{equation}
\end{theorem}

\begin{rem}\label{rem4}
The arguments $\mathbf{s}^{'}(T)W_{c}^{-1}(T)\mathbf{s}(T)$ and $\mathbf{r}^{'}(T)W_{c}^{-1}(T)\mathbf{r}(T)$ in (\ref{eq8}) are the minimum energy required to drive the system (\ref{eq1}) from $\mathbf{x}_{r}$ to $\mathbf{x}_{0}$ and $\mathbf{x}_{0}$ to $\mathbf{x}_{f}$ respectively \cite{KS72,CSL14}.
\end{rem}

From the above remark, the intuition behind Theorem \ref{thm1} is seemingly straightforward -- for a minimally novel solution to exist, the input energy must at least exceed the minimum energy required to effect the state transfer.  However, in establishing the closed-form solution for the minimally novel input, it is useful to complete a formal proof.

\begin{IEEEproof}
Define $y(t)$ as  
\begin{equation}
	y(t) = \frac{1}{T}\int_{0}^{t}\| \mathbf{u}(\tau)\|_{2}^{2}\rm{d}\tau 
\label{eq11}
\end{equation}
Clearly, $y(0) = 0$ and $y(T) = \gamma_{u}$ from (\ref{eq6b}). Thus, we can replace the constraint (\ref{eq6b}) by
\begin{equation}
	y(T) = \gamma_{u}
\label{eq12}
\end{equation}
\noindent In differential form, we can write (\ref{eq11}) as
\begin{equation}
	\frac{\rm{d}y(t)}{\rm{d}t} = \frac{1}{T}\| \mathbf{u}(t)\|_{2}^{2}
\label{eq13} 
\end{equation}

\noindent To solve the dynamic optimization problem (\ref{eq6a}), (\ref{eq6c}), and (\ref{eq12}) in continuous time, we write the Hamiltonian $\mathcal{H}(\mathbf{x}(t),y(t),\mathbf{u}(t),\mathbf{\lambda}(t),\mu(t),t)$ as
\begin{multline}
	\mathcal{H}(\mathbf{x}(t),y(t),\mathbf{u}(t),\mathbf{\lambda}(t),\mu(t),t) = -\frac{1}{T\sqrt{\gamma_{v}\gamma_{u}}} \mathbf{v}^{'}(t-T)\mathbf{u}(t) \\
	+ \mathbf{\lambda}^{'}(t)(A(t)\mathbf{x}(t)+B(t)\mathbf{u}(t)) \\
	+ \frac{\mu(t)}{T}\| \mathbf{u}(t)\|_{2}^{2}
\label{eq14}
\end{multline}
\noindent Here, $\mathbf{\lambda}(t)$ and $\mu(t)$ are the costate variables associated with the dynamics (\ref{eq1}) and (\ref{eq13}) respectively. We derive the following optimality conditions (i.e. the Euler-Lagrange equations \cite{K04}): 
\begin{subequations}
	\begin{equation}
		\frac{\rm{d}\mathbf{\lambda}(t)}{\rm{d}t} = -(\frac{\partial \mathcal{H}(\mathbf{x}(t),y(t),\mathbf{u}(t),\mathbf{\lambda}(t),\mu(t),t)}{\partial \mathbf{x}(t)})^{'} = -A^{'}(t)\mathbf{\lambda}(t)
	\label{eq15a}
	\end{equation}
	\begin{equation}
		\frac{\rm{d} \mu(t)}{\rm{d}t} = -\frac{\partial  \mathcal{H}(\mathbf{x}(t),y(t),\mathbf{u}(t),\mathbf{\lambda}(t),\mu(t),t)}{\partial y(t)}= 0
	\label{eq15b}
	\end{equation}
	\begin{equation}
		\begin{split}
			\frac{\partial \mathcal{H}(\mathbf{x}(t),y(t),\mathbf{u}(t),\mathbf{\lambda}(t),\mu(t),t)}{\partial \mathbf{u}(t)} &= 0 \\
		&= \frac{2\mu(t)}{T}\mathbf{u}(t)\\
		& - \frac{1}{T\sqrt{\gamma_{v}\gamma_{u}}}\mathbf{v}(t-T) \\
		&+ B^{'}(t)\lambda(t)
	\end{split}
	\label{eq15c}
	\end{equation}
\label{eq15}
\end{subequations}

\noindent By integrating the costate equations (\ref{eq15a}) and (\ref{eq15b}) over $t$, we obtain
\begin{subequations}
	\begin{equation}
		\mathbf{\lambda}(t) = \Phi^{'}(0,t) \mathbf{\lambda}(0)
	\label{eq16a}
	\end{equation}
	\begin{equation}
		\mu(t) \equiv \mu \qquad \forall t\in[0,T]
	\label{eq16b}
	\end{equation}
\label{eq16}
\end{subequations}

\noindent Here, $\mathbf{\lambda}(0)$ is the initial condition (at $t=0$) of (\ref{eq15a}). From (\ref{eq15c}), (\ref{eq16a}) and (\ref{eq16b}), we derive the optimal control law as
\begin{equation}
	\mathbf{u}(t) = \frac{1}{2\mu\sqrt{\gamma_{v}\gamma_{u}}}\mathbf{v}(t-T) - \frac{T}{2\mu}B^{'}(t)\Phi^{'}(0,t)\mathbf{\lambda}(0)
\label{eq17}
\end{equation}

\noindent By substituting (\ref{eq17}) into (\ref{eq6c}), we obtain $\mathbf{\lambda}(0)$ as
\begin{equation}
	\mathbf{\lambda}(0) = \frac{2\mu}{T}\Phi^{'}(T,0)W_{c}^{-1}(T)(\frac{1}{2\mu\sqrt{\gamma_{v}\gamma_{u}}}\mathbf{s}(T) - \mathbf{r}(t))
\label{eq18}
\end{equation}

\noindent  By substituting (\ref{eq17}) and (\ref{eq18}) into (\ref{eq12}), we obtain
\begin{equation}
	\mu = \pm\frac{1}{2\sqrt{\gamma_{v}\gamma_{u}}} \sqrt{\frac{\gamma_{v}T  - \mathbf{s}^{'}(T)W_{c}^{-1}(T)\mathbf{s}(T)}{\gamma_{u}T -\mathbf{r}^{'}(T)W_{c}^{-1}(T)\mathbf{r}(T)}}
\label{eq19}
\end{equation}

\noindent For the existence of a solution, $\mu$ must be a real number. Thus, either $T < \min{\{\frac{\mathbf{s}^{'}(T)W_{c}^{-1}(T)\mathbf{s}(T)}{\gamma_{v}},\frac{\mathbf{r}^{'}(T)W_{c}^{-1}(T)\mathbf{r}(T)}{\gamma_{u}}\}}$ or $T > \max{\{\frac{\mathbf{s}^{'}(T)W_{c}^{-1}(T)\mathbf{s}(T)}{\gamma_{v}},\frac{\mathbf{r}^{'}(T)W_{c}^{-1}(T)\mathbf{r}(T)}{\gamma_{u}}\}}$. Now it follows directly from Remark \ref{rem4} that the total energy $T$ must satisfy (\ref{eq8}) for the existence of a solution i.e. $T > \max{\{\frac{\mathbf{s}^{'}(T)W_{c}^{-1}(T)\mathbf{s}(T)}{\gamma_{v}},\frac{\mathbf{r}^{'}(T)W_{c}^{-1}(T)\mathbf{r}(T)}{\gamma_{u}}\}}$. 

\end{IEEEproof}

\subsubsection{Uniqueness of the Minimally Novel Input}
\begin{theorem}\label{thm2}
	Under the hypothesis of Theorem \ref{thm1}, the solution of the non-convex optimization problem (\ref{eq6}) is unique.
\end{theorem}

\begin{IEEEproof}
\noindent By replacing $\lambda(0)$ in (\ref{eq17}) with (\ref{eq18}), then substituting into (\ref{eq3}), we obtain the optimal value of $\mathbb{J}(T)$ as a function of $\mu$:
\begin{multline}
	\mathbb{J}(T) = \frac{1}{T\sqrt{\gamma_{v}\gamma_{u}}}\mathbf{s}^{'}(T)W_{c}^{-1}(T)\mathbf{r}(T)\\
	+\frac{1}{2\mu\gamma_{u}}(1-\frac{1}{\gamma_{v}T}\mathbf{s}^{'}(T)W_{c}^{-1}(T)\mathbf{s}(T))
\label{eq20}
\end{multline}

\noindent It follows from Theorem \ref{thm1} that $\frac{1}{\gamma_{v}T}\mathbf{s}^{'}(T)W_{c}^{-1}(T)\mathbf{s}(T) \in (0,1)$ (see (\ref{eq8})). Thus, the maximum of $\mathbb{J}(T)$ occurs when $\mu > 0$ in (\ref{eq19}) i.e. 
\begin{equation}
	\mu = \frac{1}{2\sqrt{\gamma_{v}\gamma_{u}}} \sqrt{\frac{\gamma_{v}T  - \mathbf{s}^{'}(T)W_{c}^{-1}(T)\mathbf{s}(T)}{\gamma_{u}T -\mathbf{r}^{'}(T)W_{c}^{-1}(T)\mathbf{r}(T)}}
\label{eq21}
\end{equation}  

\noindent Thus, a unique optimal control input $\mathbf{u}(t)$ exists and is given by
\begin{multline}
	\mathbf{u}(t) = \frac{1}{2\mu\sqrt{\gamma_{v}\gamma_{u}}}(\mathbf{v}(t-T)-B^{'}(t)\Phi^{'}(T,t)W_{c}^{-1}(T)\mathbf{s}(T)) \\
	+ B^{'}(t)\Phi^{'}(T,t)W_{c}^{-1}(T)\mathbf{r}(T)
\label{eq22}
\end{multline}
\end{IEEEproof}

The expression \eqref{eq22} is the closed form for the minimally novel input.  We now proceed to state an important lemma that establishes the dependence of the minimum novelty control on only the initial, intermediate and terminal states.  
\begin{lemma}\label{lem1}
	Given $\mathbf{x}_{r}$,  $\mathbf{x}_{0}$ and  $\mathbf{x}_{f}$, $\mathbb{J}(T)$ in (\ref{eq20}) is independent of $\mathbf{v}(t-T)$, $t\in[-T,0)$. 
\end{lemma}
\begin{IEEEproof}
Since the input $\mathbf{v}(t-T)$ drives (\ref{eq1}) from $\mathbf{x}_{r}$ to $\mathbf{x}_{0}$ in time $T$, we can rewrite (\ref{eq9a}) as
\begin{equation}
	\mathbf{s}(T) = \int_{0}^{T}\Phi(T,t)B(t)\mathbf{v}(t-T)\rm{d}t = \mathbf{x}_{0}-\Phi(T,0)\mathbf{x}_{r}
\end{equation}
Clearly, $\mathbf{s}(T)$ is independent of $\mathbf{v}(t-T)$. Thus, $\mathbb{J}(T)$ in (\ref{eq20}) is independent of $\mathbf{v}(t-T)$.
\end{IEEEproof}
This lemma will form the basis for our use of novelty as a control metric, as highlighted in Section \ref{examples} below.

\subsubsection{Euclidean - Inner Product Transformation}
\noindent As noted in Remark \ref{rem3}, it is a notable consequence of our cost formulation that the problem can be exactly recast in terms of the Euclidean norm. Specifically, if we consider 

\begin{equation}
\mathbb{J}_{1}(T) = \frac{1}{T}\int_{0}^{T} \|\mathbf{v}(t-T)-\mathbf{u}(t)\|_{2}^{2}\rm{d}t
\label{eq23}
\end{equation}

\noindent  as the cost function in (\ref{eq6a}), we obtain the optimal solution as
\begin{subequations}
	\begin{equation}
		\mu = -1 + \sqrt{\frac{\gamma_{v}T  - \mathbf{s}^{'}(T)W_{c}^{-1}(T)\mathbf{s}(T)}{\gamma_{u}T -\mathbf{r}^{'}(T)W_{c}^{-1}(T)\mathbf{r}(T)}}
	\label{eq24a}
	\end{equation}
	\begin{multline}
		\mathbf{u}(t) = \frac{1}{1+\mu}(\mathbf{v}(t-T)-B^{'}(t)\Phi^{'}(T,t)W_{c}^{-1}(T)\mathbf{s}(T))\\
	+ B^{'}(t)\Phi^{'}(T,t)W_{c}^{-1}(T)\mathbf{r}(T)
	\label{eq24b}
	\end{multline}
	\begin{multline}
		\mathbb{J}_{1}(T) = (\gamma_{u}+\gamma_{v}) - \frac{2}{T}\mathbf{s}^{'}(T)W_{c}^{-1}(T)\mathbf{r}(T))\\
		+\frac{2\gamma_{v}}{1+\mu}(\frac{1}{\gamma_{v}T}\mathbf{s}^{'}(T)W_{c}^{-1}(T)\mathbf{s}(T)-1)
	\label{eq24c}
	\end{multline}
\label{eq24}
\end{subequations} 

\noindent It is evident that the control law (\ref{eq21})-(\ref{eq22}) is same as the control law (\ref{eq24a})-(\ref{eq24b}), as one expects from Remark \ref{rem3}.   

\subsection{Average Novelty Measure and Results}
\noindent The definition of the inner-product novelty measure (\ref{eq3}) requires the inputs $\mathbf{v}(t-T)$ and $\mathbf{u}(t)$ to be of equal time length $T$.  To allow for characterization of the input novelty in cases where the time lengths of these inputs are different, we define an average novelty measure in the form  

\begin{equation}
\mathbb{J}_{2}(T) = (\frac{1}{\sqrt{\gamma_{v}}T^{*}}\int_{0}^{T^{*}} \mathbf{v}(t-T^{*})\rm{d}t)^{'}(\frac{1}{\sqrt{\gamma_{u}}T}\int_{0}^{T} \mathbf{u}(t)\rm{d}t)
\label{eq25}
\end{equation}

\noindent with 

\begin{equation}
\frac{1}{T^{*}}\int_{0}^{T^{*}}\|\mathbf{v}(t-T^{*})\|_{2}^{2}\rm{d}t = \gamma_{v} 
\label{eq26}
\end{equation}

\noindent Here, $T^{*}$ is the time duration of the prior input $\mathbf{v}(t-T^{*})$ which drives the system (\ref{eq1}) from $\mathbf{x}_{r}$ to $\mathbf{x}_{0}$. If we consider (\ref{eq25}) as the cost function in (\ref{eq6a}), we obtain the optimal $\mathbf{u}(t)$ as 

\begin{multline}
	\mathbf{u}(t) =\frac{1}{2\mu\sqrt{\gamma_{v}\gamma_{u}}}(\mathbf{v}_{av}-B^{'}(t)\Phi^{'}(T,t)W_{c}^{-1}(T)\mathbf{s}(T))\\
	+ B^{'}(t)\Phi^{'}(T,t)W_{c}^{-1}(T)\mathbf{r}(T)
\label{eq27}
\end{multline}
\noindent where
\begin{subequations}
	\begin{equation}
		\mathbf{s}(T) = \int_{0}^{T}\Phi(T,t)B(t)\mathbf{v}_{av}\rm{d}t
	\label{eq28a}
	\end{equation}
	\begin{equation}
		\mathbf{v}_{av} = \frac{1}{T^{*}}\int_{0}^{T^{*}} \mathbf{v}(t-T^{*})\rm{d}t
	\label{eq28b}
	\end{equation}
\label{eq28}
\end{subequations}

\noindent The optimal $\mathbb{J}_{2}(T)$ and $\mu$ are given by (\ref{eq20}) and (\ref{eq21}) respectively. It should be noted here that, in this case, $\mathbf{s}^{'}(T)W_{c}^{-1}(T)\mathbf{s}(T)$ is no longer the minimum energy required to drive the system (1) from $\mathbf{x}_{r}$ to $\mathbf{x}_{0}$ (see Remark \ref{rem4}).

\section{Discrete-Time, Linear Dynamical Networks}\label{sec4}
\subsection{Input Novelty}
\noindent We consider linear, time-varying discrete time systems with dynamics of the form
\begin{equation}
	x(k+1) = A(k)x(k) + B(k)u(k)
\label{eq29}
\end{equation} 
\noindent Here $x(k)\in \mathbb{R}^{n\times 1}$ represents the state of the system at discrete time $k$ for $k = 0,1,2,\ldots$, $A(k)\in \mathbb{R}^{n\times n}$ is the time-varying dynamic matrix, $B(k)\in \mathbb{R}^{n\times m}$ is the input matrix, and $\mathbf{u}(k)\in \mathbb{R}^{m\times 1}$ is the input to the system. We say that (\ref{eq29}) describes the time evolution of discrete time linear networks in the presence of external inputs.    

Let us assume an input $\mathbf{v}(k-p) \in \mathbb{R}^{m\times 1}$, $k=0,1,2,\ldots,p-1$, with total energy $p\gamma_{v}$, i.e.
\begin{equation}
	\frac{1}{p}\sum_{k=0}^{p-1}\|\mathbf{v}(k-p)\|_{2}^{2} = \gamma_{v}
\label{eq30}
\end{equation} 

\noindent Here, $\gamma_{v}>0$ is a constant scaling factor. We assume that the sequence of inputs $\mathbf{v}(k-p), k=0,1,2,\ldots,p-1$ can drive $\mathbf{x}(k)$ from $\mathbf{x}_{r}$ to $\mathbf{x}_{0}$, subject to the dynamics (\ref{eq29}). We introduce the inner product   

\begin{equation}
	\mathbb{J}(p) = \frac{1}{p\sqrt{\gamma_{u}\gamma_{v}}}\sum_{k=0}^{p-1} \mathbf{v}^{'}(k-p)\mathbf{u}(k)
\label{eq31}
\end{equation}

\noindent where
\begin{equation}
	\frac{1}{p}\sum_{k=1}^{p}\|\mathbf{u}(k)\|_{2}^{2} = \gamma_{u}
\label{eq32}
\end{equation}

\noindent to measure the novelty of a subsequent input $\mathbf{u}(k)$, $k=0,1,2,\ldots,p-1$, relative to $\mathbf{v}(k-p)$, required in order to reach the state $\mathbf{x}_{f}$. In other words for a fixed input energy, $\mathbb{J}(p)$ measures the required directional change in inputs (thus novelty) to achieve a given state change in the state of the system. Here, $\gamma_{u}>0$ is a constant scaling factor. 

\begin{rem}\label{rem5}
	Remarks \ref{rem1}, \ref{rem2} and \ref{rem3} are applicable for discrete time systems if we replace $t$ and $T$ by $k$ and $p$ respectively. 
\end{rem}

\subsection{Minimum Novelty Problem}
\noindent Similar to the continuous-time case, we formulate the following optimal control problem to design the minimally novel input $\mathbf{u}(k), k=0,1,2,\ldots,p-1$ such that a desired change in the state of the system can be achieved under the constraint of fixed energy subject to the system dynamics (\ref{eq29}):

\begin{subequations}
\begin{equation}
    \min_{\substack{\mathbf{u}(k) \\k=1,2,\ldots,p}}
         \quad   -\mathbb{J}(p) \label{eq33a} 
\end{equation}
    \textrm{s.t.} 
\begin{equation}
        \frac{1}{p} \sum_{k=0}^{p-1}\| \mathbf{u}(k)\|_{2}^{2}= \gamma_{u} \label{eq33b}
\end{equation}
\begin{multline}
	\mathbf{x}_{f} = \prod_{k=0}^{p-1}A(p-k-1)\mathbf{x}_{0} \\
+ \sum_{k=0}^{p-1}\prod_{j=k+1}^{p-1}A(p-j+k)B(k)\mathbf{u}(k)  \label{eq33c}
\end{multline}
\label{eq33}
\end{subequations}

\noindent It should be noted here that the constraint (\ref{eq33c}) is obtained by computing $\mathbf{x}_{f}$ using the dynamics (\ref{eq29}). Immediately, we note that the quadratic equality constraint (\ref{eq33b}) makes the optimization problem (\ref{eq33}) non-convex.

\subsection{Results}
\noindent 
We relax the non-convex problem by replacing the equality constraint \eqref{eq33b} with an inequality constraint.  Thus, the relaxed problem becomes convex and by deriving the optimality conditions, we establish that the relaxation is exact, i.e. any solution of the non-convex optimization problem (\ref{eq33}) is also the solution of the convex optimization problem.  

\subsubsection{Convex Relaxation}
\noindent We relax the equality constraint (\ref{eq33b}) with the following inequality constraint 

\begin{equation}
	\frac{1}{p} \sum_{k=0}^{p-1}\| \mathbf{u}(k)\|_{2}^{2} \leq \gamma_{u} 
\label{eq34}
\end{equation}

\noindent It is clear that the optimization problem (\ref{eq33a}), (\ref{eq33c}), and (\ref{eq34}) is convex. 

\subsubsection{Exact Relaxation}
\begin{lemma}
	A solution of the non-convex optimization problem (\ref{eq33}) is also the solution of the convex optimization problem (\ref{eq33a}), (\ref{eq33c}), and (\ref{eq34}). 
\end{lemma}

\begin{IEEEproof}
	It is sufficient to show that equality holds for the constraint (\ref{eq34}) in the convex optimization problem (\ref{eq33a}), (\ref{eq33c}), and (\ref{eq34}). For this, we write the Lagrangian as
\begin{multline}
	\mathcal{L}(\gamma,\mathbf{\delta},\mathbf{u}(k), k=0,1,\ldots,p-1) = -\mathbb{J}(p) \\
+ \gamma( \frac{1}{p} \sum_{k=0}^{p-1}\| \mathbf{u}(k)\|_{2}^{2}-\gamma_{u}) + \sum_{i=1}^{n}\delta_{i}(\sum_{l=1}^{n}(\prod_{k=0}^{p-1}A(p-k-1))_{il}(\mathbf{x}_{0})_{l}) \\ 
+ \sum_{i=1}^{n}\delta_{i}(\sum_{l=1}^{n}(\sum_{k=0}^{p-1}\prod_{j=k+1}^{p-1}A(p-j+k)B(k))_{il}u_{l}(k)-(\mathbf{x}_{f})_{i})
\label{eq35}
\end{multline}

\noindent Here, $\lambda$ and $\mathbf{\delta}$ are the Lagrange multipliers associated with the constraints (\ref{eq34}) and (\ref{eq33c}) respectively. We derive the Karush-Kuhn-Tucker (KKT) conditions for optimality as:

\begin{subequations}
	\begin{multline}
		\frac{\partial \mathcal{L}(\gamma,\mathbf{\delta},\mathbf{u}(k), k=0,1,\ldots,p-1)}{\partial u_{j}(k)} = 0 \qquad j = 1,2,\ldots,m \\
		\qquad k = 0,1,\ldots,p-1
	\label{eq36a}
	\end{multline}
	\begin{equation}
		 \frac{1}{p} \sum_{k=0}^{p-1}\sum_{j=1}^{m}u_{j}^{2}(k)\leq \gamma_{u}
	\label{eq36b}
	\end{equation}
	\begin{equation}
		\gamma \geq 0
	\label{eq36c}
	\end{equation}
	\begin{equation}
		\gamma( \frac{1}{p} \sum_{k=0}^{p-1}\sum_{j=1}^{m}u_{j}^{2}(k)-\gamma_{u}) = 0
	\label{eq36d}
	\end{equation}
	\begin{multline}
		\mathbf{x}_{f} = \prod_{k=0}^{p-1}A(p-k-1)\mathbf{x}_{0} \\
+ \sum_{k=0}^{p-1}\prod_{j=k+1}^{p-1}A(p-j+k)B(k)\mathbf{u}(k)  
	\label{eq36e}
	\end{multline}
\label{eq36}
\end{subequations}

\noindent From (\ref{eq35}) and (\ref{eq36a}), we obtain the optimal control law as
\begin{equation}
	u_{j}(k) = \frac{\frac{1}{\sqrt{\gamma_{u}\gamma_{v}}}v_{j}(k-p) - p\sum_{i=1}^{n}\delta_{i}\prod_{l=k+1}^{p-1}A(p-l+k)B(k))_{ij}}{2\gamma}
\label{eq37}
\end{equation}

\noindent It follows from (\ref{eq36b}), (\ref{eq36c}), and (\ref{eq37}) that a solution of the optimization problem exists if $\gamma > 0$. Using this along with (\ref{eq36d}), we find that 
\begin{equation}
	\frac{1}{p} \sum_{k=1}^{p}\sum_{j=1}^{m}u_{j}^{2}(k) = \gamma_{u}
\label{eq38}
\end{equation}
\end{IEEEproof}
The result permits numerically efficient evaluation of the minimally novel input.

\section{Minimum Novelty as a Control Metric}\label{examples}
\noindent We now return to the original motivation outlined in the introduction: to systematically characterize the responsiveness of linear networks to novelty.  For instance, given systems $S_1$ and $S_2$, identical in structure but different in dynamics, we may ask which requires more novelty in order to execute the same trajectory.  In this sense, we can fashion the minimum novelty solutions analytically obtained above into a control (or, controllability) metric.  To do so, we standardize the notions of initial, intermediate and terminal states.  The space of terminal states can then be sampled systematically (e.g., along eigenvectors of the Controllability gramian), or via some other Monte-Carlo fashion, in order to gain an aggregate novelty-based metric.

\subsection{Examples}\label{sec5}

\begin{figure}
	\centering
		\includegraphics[width=0.48\textwidth]{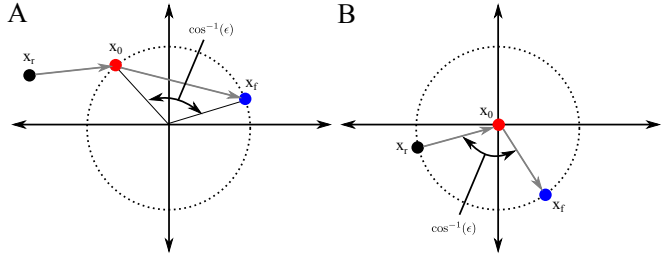}
	\caption{Schematic of initial, intermediate and final states for (A) Example 1, (B) Example 2.  The dashed circle has unit radius.}
	\label{fig:example_setup}
\end{figure}

\subsubsection{Rate-based Recurrent Neuronal Networks}\label{sec5a}
\noindent We consider a recurrent network of $n$ neurons with linearized firing rate dynamics of the form \cite{DA01} 
\begin{equation}
	S\frac{\rm{d}\mathbf{x}(t)}{\rm{d}t} = -\mathbf{x}(t) + W \mathbf{x}(t) + B\mathbf{u}(t)
\label{eq39}
\end{equation}

\noindent Here, $\mathbf{x}(t) \in \mathbb{R}_{+}^{n\times 1}$ represents the firing rate of the neurons at time $t$, $S \in \mathbb{R}_{+}^{n\times n}$ is a diagonal matrix whose diagonal elements are the (positive) time constants of the neurons, $W \in \mathbb{R}^{n\times n}$ defines the linear interaction among neurons in the network, $B \in \mathbb{R}_{+}^{n\times n}$ is the input matrix, and $\mathbf{u}(t)\in \mathbb{R}_{+}^{n\times 1}$ is the afferent input. Since $S$ is invertible, (\ref{eq39}) can be represented in the form of (\ref{eq1}) by considering $A = S^{-1}(-I+W)$ where $I$ is the $n\times n$ identity matrix. 

In this example, we consider a recurrent network of $n = 100$ neurons where $80$ neurons are excitatory and $20$ are inhibitory. We choose the time constants (in milliseconds) of the neurons, i.e. the elements of the diagonal matrix $S$, from a uniform distribution $\mathcal{U}(5,10)$. 
The connectivity weights $w_{i,j}$ (in essence, a time constant for excitation/inhibition from the neuron $i$ to $j$) are drawn from a uniform distribution $\mathcal{U}(0,1)$ or $\mathcal{U}(-1,0)$, if neuron $i$ is excitatory or inhibitory, respectively.
We assume that $w_{i,j} = 0$ for $i=j$, i.e. neurons do not directly excite/inhibit themselves.
Assuming $B = S$, we proceed to compute the minimally novel inputs required to make a desired directional change in the network firing rate using (\ref{eq8})-(\ref{eq10}), (\ref{eq20})-(\ref{eq22}). 

To complete the example, we specify $T = 3$ ms. The states $\mathbf{x}_{0}$ and $\mathbf{x}_{f}$ are specified to satisfy $\|\mathbf{x}_{0}\|_{2} = \|\mathbf{x}_{f}\|_{2} = 1$ with $\mathbf{x}_{0}'\mathbf{x}_{f} = \epsilon$, where in this particular case we specify $\epsilon = 0.7645$. The average energy $\gamma_{v}$ of the prior input $\mathbf{v}(t-T)$ (see (\ref{eq2})) and the average energy $\gamma_{u}$ of the computed minimally novel input of $\mathbf{u}(t)$ (see (\ref{eq4})) are set to $1$. The prior input $\mathbf{v}(t-T)$ is specified to be constant over the interval $t\in[0,T]$.  Figure \ref{fig:example_setup}A illustrates the schematic of the example setup.

\begin{figure}
\centering%
\includegraphics[width=0.48\textwidth]{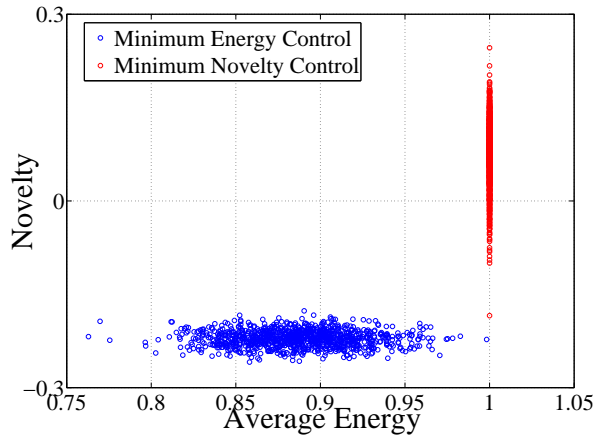}
\caption{Comparison of minimum novelty control with minimum energy control for $n=1000$ random realizations of the recurrent neuronal network: Each red dot on the figure depicts the novelty associated with the solution to the minimum novelty control. Each blue dot corresponds to the minimum energy solution.}
\label{fig3}
\end{figure}

Figure \ref{fig3} illustrates the outcome of the example for $n=1000$ random realizations of the system. Each red dot on the figure depicts the novelty associated with the solution to (\ref{eq20})-(\ref{eq21}), i.e., the minimum novelty. Note, again, that in this example, these inputs are constrained to have unit average energy. Each blue dot corresponds to the minimum energy solution (for the same $\mathbf{x}_{0}$ and $\mathbf{x}_{f}$). As a verification of our theoretical development, we note that the minimum energy solution consistently requires an injection of novelty (angular orientation) relative to the prior input and relative to the optimum.

\subsection{Complex Networks}\label{sec5b}   
\noindent For our second example, we demonstrate the utility of our novelty-based controllability metric for distinguishing linear complex networks based on topological properties. Moreover, we demonstrate that a novelty-based control metric provides a different description of network controllability than does an energy-based characterization.   

Consider networks of $n=100$ excitatory neurons with dynamics of the form (\ref{eq39}) where the structure of the weight matrix $W$ is defined by the adjacency matrix of either the undirected Barab\'{a}si-Albert scale-free network (BA) \cite{BA99} or the undirected Watts-Strogatz small-world network (WS) \cite{WS98} with appropriate parameters. If the $(i,j)^{th}$ element of the network adjacency matrix is $1$, we say that neuron $i$ is connected to the neuron $j$. 

For both types of networks, we construct the matrices $S$ and $W$ (see (\ref{eq39})) as follows. Given a degree distribution of the BA network, we construct $1000$ full rank realizations of the adjacency matrix. It should be noted here that the total number of undirected edges in all $1000$ realizations is the same. For each realization, we replace the nonzero elements of the adjacency matrix by a number drawn from a uniform distribution $\mathcal{U}(0,1)$, thus yielding the matrix $W$.  We choose the time constants (in milliseconds) of the neurons (the elements of the diagonal matrix $S$) from a uniform distribution $\mathcal{U}(5,10)$. An analagous procedure is followed to specify the WS networks. 
Since the development of a WS network does not involve the creation of new edges, but rather the `rewiring' of existing connections, while a BA network is created by a growth algorithm, we adjust their respective parameters to obtain approximately the same number of edges in both networks.

Assuming $B = S$, we proceed to compute the minimum directional change in inputs (i.e. minimally novel inputs) required to make a desired directional change in firing rates of neurons using (\ref{eq8})-(\ref{eq9}), (\ref{eq19})-(\ref{eq21}). 

To complete the example, we again specify $T = 0.3$ ms. The states $\mathbf{x}_{r}$ and $\mathbf{x}_{f}$ are specified to satisfy $\|\mathbf{x}_{r}\|_{2} = \|\mathbf{x}_{f}\|_{2} = 1$ with $\mathbf{x}_{r}'\mathbf{x}_{f} = \epsilon$, where in this particular case we specify $\epsilon = 0.7358$. We set the intermediate state $\mathbf{x}_{0}$ to $\mathbf{0}$. As a reminder, the prior input $\mathbf{v}(t-T)$ drives the network from $\mathbf{x}_{r}$ to $\mathbf{x}_{0}$ in time $T$ and the computed minimally novel input of $\mathbf{u}(t)$ drives the network from $\mathbf{x}_{0}$ to $\mathbf{x}_{f}$ in time $T$. The average energy $\gamma_{v}$ of $\mathbf{v}(t-T)$ (see (\ref{eq2})) and the average energy $\gamma_{u}$ of $\mathbf{u}(t)$ (see (\ref{eq4})) are set to $300$.  Figure \ref{fig:example_setup}B illustrates the example setup.  Note that from Lemma \ref{lem1}, it is sufficient to randomly sample the system parameters without explicitly sampling the space of prior inputs.

\begin{figure}
\centering%
\includegraphics[width=0.48\textwidth]{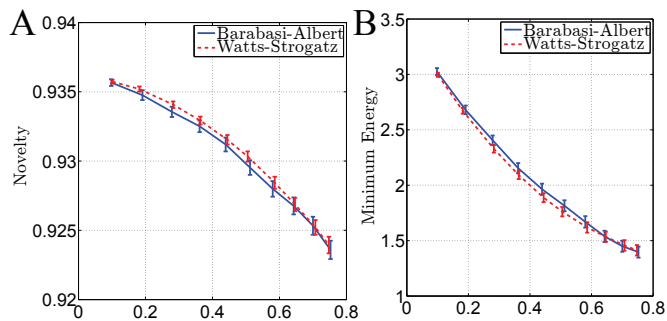}
\caption{Controllability based distinction between the Barabasi-Albert scale free network (blue) and the Watts-Strogatz small world network (red). (A) shows the novelty associated with the solution to the minimum novelty control as the function of the total number of edges in the network. (B) shows the required minimum energy as the function of the total number of edges in the network. For each relevant $\eta$, we show the mean and mean $\pm$ standard deviation for $1000$ realizations of the network.}
\label{fig4}
\end{figure}

Figure \ref{fig4}A shows the minimum novelty associated with the solution of (\ref{eq20})-(\ref{eq21}) as a function of the total number of edges in the BA network (blue) and the WS network (red). Here, $\eta$ is defined as the ratio of the total number of undirected edges in the network and the total possible number of undirected edges (i.e. $\frac{n(n-1)}{2}$). At each $\eta$, we show the mean and mean $\pm$ standard deviation over $1000$ realizations of the network. As shown in this figure, (i) controlling both types of networks becomes progressively harder in terms of minimum novelty (more novelty is required) as these networks become more dense, and (ii) the WS network requires consistently less novelty to control.

Figure \ref{fig4}B shows the minimum energy required to effect the same state trajectories for the BA network (blue) and the WS network (red). As shown in this figure, (i) controlling both types of networks becomes progressively easier in terms of minimum energy (less energy is required) as these networks become more dense, and (ii) the WS network requires consistently less energy to control. 

This example demonstrates how a novelty-based analysis can provide a different and important controllability characterization, complementing the use of minimum energy-based methods. Moreover, both metrics are consistent in distinguishing the BA network from the WS network.

\section{Conclusions}
\noindent In many natural systems, and especially in biological sensory networks, the responsiveness of the system to a given input is dependent on both the energy of the stimulus and also its salience or novelty relative to background and prior inputs.  Conventional control analysis aggregates all past input history into the current state and evaluates controllability in terms of the forward trajectory.  Here, we introduce a control analysis to explicitly quantify the extent of novelty, defined in terms of an inner product, required to effect a state transfer in a linear system.  The analysis, and corresponding solution to the problem of finding minimally novel inputs, are used to establish a new type of characterization of system controllability.  Through two illustrative examples, we show how such an analysis is distinct from energy-based metrics in describing the controllability of large-scale linear systems or networks.  The results indicate that the analysis of such systems may benefit from a multidimensional metric in which energy, novelty and other control characteristics are components.  In future work, we plan to use such an approach in the specific domain of sensory neuroscience to examine how various time-scales of dynamics, including neuroplasticity, impact the ability of networks to transform afferent input and facilitate higher information processing.  Application in other domains, including finance and social network analysis, could similarly be conceived.






%



\section*{Acknowledgment}
\noindent S. Ching Holds a Career Award at the Scientific Interface from the Burroughs-Wellcome Fund.

\ifCLASSOPTIONcaptionsoff
  \newpage
\fi



\bibliographystyle{IEEEtran}
\bibliography{masterN}
\end{document}